\documentclass{amsart}
\usepackage{amsfonts,verbatim,epigraph}
\usepackage[arrow,matrix,curve,ps]{xy}

\begin{document}

\title{Richard Dedekind: style and influence}
\author{Stefan M\"uller-Stach}
\address{Institut f\"ur Mathematik, Staudingerweg 9, Johannes Gutenberg--Universit\"at Mainz, 55099 Mainz, Germany, email: mueller-stach@uni-mainz.de}


\begin{abstract} 
This text is based on an invited talk at the Dedekind memorial conference at Braunschweig in October 2016. 
It summarizes views from \cite{sms}. 
\end{abstract}

\maketitle

\section{Intro} 

The monograph \cite{sms} contains a commented edition of the two books ``Stetigkeit und Irrationale Zahlen'' \cite{dedekind1872} (1872)  
and ``Was sind und was sollen die Zahlen?'' \cite{dedekind1888} (1888) by Richard Dedekind. 
Both books were conceived much earlier, the first in the late 1850's as a preparation for lectures in Z\"urich 
and the second in the 1870's, as handwritten manuscripts in the G\"ottingen archives show. 

Hilbert reported -- after discussing with other mathematicians -- about 
the hostile (``gegnerische'') reception of Dedekind's 1888 book \cite[Vorwort]{sms}: 

\begin{quote}
``Im Jahre 1888 machte ich als junger Privatdozent von K\"onigsberg aus eine Rundreise an die deutschen Universit\"aten. Auf meiner ersten 
Station, in Berlin, h\"orte ich in allen mathematischen Kreisen bei jung und alt von der damals erschienenen Arbeit Dedekinds ,,Was sind und was sollen die Zahlen?'' 
sprechen -- meist in gegnerischem Sinne. Die Abhandlung ist neben der Untersuchung von Frege der wichtigste erste tiefgreifende Versuch einer Begr\"undung der elementaren Zahlenlehre.'' 
\end{quote}

At the time, mathematicians had similar problems while digesting Dedekind's equally influential and innovative work on Galois theory and algebraic number theory. 
Dedekind was so much ahead of his time in his research that some of his achievements reach without any 
necessary ``historical transformation'' into the 20th century. The editors of Dedekind's collected works (Fricke, Noether, Ore) confirmed this view in 1932 \cite[Nachwort, Vol. 3]{dedekind}: 

\begin{quote}
``Es ist ein Zeichen, wie Dedekind seiner Zeit voraus war, da{\ss} seine Werke noch heute lebendig sind, ja da{\ss} sie vielleicht erst heute ganz lebendig geworden sind.''
\end{quote}

Dedekind's style has many implications for the reception of his work, as we will try to explain here. 

\section{Character and style}

Dedekind had a very close relationship with his whole family. He lived with his mother and his sister Julie in their family's home until Julie's death in 1914.
We imagine the siblings spending hours together reading and talking. The family owned a weekend retreat in Harzburg where they used to meet 
with other friends and went for hiking tours. Richard was close to the family of his brother Adolf as well and supported them during a period of health problems of his nephew Atta. 
These strong ties made Dedekind never leave his hometown Braunschweig for other universities after he had returned from his first position in Z\"urich.

Dedekind had an impressive and unique style in his works and showed a generous character in his communication with colleagues. 
He characterized his own mathematical style in the following way \cite[Preface]{dedekind1888}: 

\begin{quote}
``Was beweisbar ist, soll in der Wissenschaft nicht ohne Beweis 
geglaubt werden ...\\
Indem ich die Arithmetik (Algebra, Analysis) nur einen Teil der Logik nenne  ... \\
Die Zahlen sind freie Sch\"opfungen des menschlichen Geistes ...'' 
\end{quote}

The second sentence shows that he considered himself being a \emph{logicist} (logicism was originated by Frege as a program to reduce mathematics to logic).  
Dedekind's words sound quite dramatic in a way. Indeed, the times in which he lived were times of major changes in science. 
Epple called the underlying phenomenon ``Das Ende der Gr\"o{\ss}enlehre'' \cite{epple}. The invention of set-theory by Cantor and Dedekind was a scientific revolution, comparable to 
the invention of Riemannian manifolds or Galois theory.

The occurence of the word \emph{Sch\"opfungen} (i.e., \emph{creations}) is remarkable. There are three areas where Dedekind employed
this technique. First, in an unpublished manuscript ``Die Sch\"opfung der Null und der negativen ganzen Zahlen'' from 1872, in which 
he constructs $\mathbb{Z}$ and $\mathbb{Q}$ from $\mathbb{N}$ using equivalence relations. Then, in the 1872 book ``Stetigkeit und 
Irrationale Zahlen'' \cite{dedekind1872}, where he invented Dedekind cuts, and finally in his creation of ``ideal numbers'' through ideals in number rings. 
In all three cases, new numbers are defined as sets. Divisibility and other concepts arise from set-theoretic operations.
Dedekind saw a step of abstraction in defining the new numbers as single objects, therefore he used the phrase ``Sch\"opfung''. 
He defended this viewpoint in a 1888 letter to H. Weber \cite[p. 276]{scheel}. 

Dedekind's mathematical style was characterized by other people as ``axiomatique, formel et abstrait'' (Sinaceur 1974) and as 
``built on traditions, abstract und Bourbaki style'' (Scharlau 1981). Similarities to the Bourbaki group and eventually to Grothendieck are quite important. 
Dedekind prefered to study mathematics in a more formal way by using mappings instead of element-wise considerations and theoretical concepts instead of examples. 
In this sense, he was quite different from mathematicians of his time.  

I would like to add that his style was very elegant, always concise and well-readable for a modern reader, so also in this respect Dedekind was far ahead of his time.
This ``modernness'' can also be expressed by saying that there is almost no ``historical transformation'' necessary in order to translate his ideas into our times. 
Dedekind published only relatively few papers after a rather long time of thinking. He would probably hate the ``modern'' rapid publication style with its enormous output.  

A programmatic document for Dedekind's approach to mathematics is his remarkable \emph{Habilitationsvortrag} from 1854 \cite[Vol. 3, p. 428]{dedekind}. 
In it he describes first some basics of scientific theories in general and then proceeds to set up his program for the creation of new objects and theories in mathematics.

Some years later, Dedekind clearly distinguishes his own views from those of others in the preface of his 1888 book \cite{dedekind1888}: 

\begin{quote}
``Das Erscheinen dieser Abhandlungen [von Helmholtz und Kronecker] ist die Veranlassung, die mich bewogen hat, nun auch mit meiner, in mancher Beziehung \"ahnlichen, aber
durch ihre Begr\"un\-dung doch wesentlich verschiedenen Auffassung hervorzutreten, die ich mir seit vielen Jahren und ohne jede Beeinflussung gebildet habe.'' 
\end{quote}

Dedekind's personal style is present in his correspondence and relationships. 
The letters to Cantor and others \cite{cavailles-noether,scheel} show -- in their friendly and open way of sharing common knowledge --  
many positive features of the personality of Dedekind. Many letters to his family were published by Ilse Dedekind \cite{ilse_dedekind}.
They prove that the older colleagues Gau{\ss}, Riemann and Lejeune Dirichlet were Dedekind's ``Vorbilder'' in G\"ottingen. 

In particular, Dedekind was close to Dirichlet, which eventually led to his edition of Dirichlet's number theory lectures with the famous appendix (the ``supplement'') 
on algebraic number theory. Dirichlet seems to have influenced Dedekind also in foundational matters in the way he prefered conceptual mathematics to examples 
(sometimes called the ``second Dirichlet principle'') and in the idea of reducing mathematics to the theory of natural numbers \cite[\S 1]{sms}. 

Dedekind had a very strong and long-lasting friendship with Heinrich Weber. Both collaborated during the 
difficult edition of Riemann's collected works, which is documented in their letters \cite{scheel}. 
In joint work they developed the theory of function fields in one variable \cite[Vol. 1, p. 238]{dedekind}. 

The relationship with Kronecker, a powerful figure in Berlin, must have been uncomfortable at least during later years.
Kronecker, in his lectures on number theory from 1891, and Weyl, in his book ``Das Kontinuum'' (1918), both
explicitly criticize Dedekind's treatment of real numbers. Dedekind appears in retrospect as one of the opponents (together with Hilbert) of 
dogmatic forms of \emph{intuitionism} and \emph{constructivism} in mathematics. 

After Cantor's unthoughtful publication of some common ideas about the countability of algebraic numbers 
without mentioning Dedekind, thus violating scientific standards boldly, 
Dedekind did not react anymore to Cantor's letters for a certain period from 1874 on, most likely because of this instance. 
A letter from Dedekind's nephew Atta to Julie from 1887 \cite[p. 223]{ilse_dedekind} proves that this was a well-known issue in the entire Dedekind family. 
In 1882 both met again, but even after this, as the 1888 letter to H. Weber \cite[p. 276]{scheel} shows, Dedekind was not happy with Cantor's attitude. 
In a letter to Hilbert from 1899 \cite[\S 1]{sms}, Cantor avoided to mention any of this and only emphasized the problems 
about the existence of an infinite set. In summary, the personalities of both were quite different and apparently sometimes incompatible.

\section{Axioms for $\mathbb{N}$ and recursion theory} 

Dedekind found appropriate ``axioms'' for the natural numbers $\mathbb{N}$ (without calling them axioms) 
and constructed models of $\mathbb{N}$ using some portion of his newly invented elementary set-theory.

Dedekind considered an infinite set $X$ together with an injective self-map (``successor map'') $S: X \to X$ and a distinguished element $0 \in X$. 
We write $(X,0,S)$ for a triple with these properties. A \emph{chain} is a subset $C \subset X$ such that $C$ is stable under $S$, i.e., $S(C) \subset C$. 
A chain containing $0$ is called \emph{simply-infinite} if it is the chain of all successors of $0$, i.e., the intersection of all chains containing $0$: 
\[
\bigcap_{C \, {\rm chain}, \,0 \in C} C =\{0,1:=S(0),2:=S(S(0)), \ldots \}.  
\] 
Simply-infinite chains are determined uniquely by the triple $(X,0,S)$. 
For Dedekind, a model $\mathbb{N}$ for the natural numbers is a simply-infinite chain arising from a triple $(X,0,S)$ where the set $X$ may be arbitrarily large. Let us 
denote the restriction of the map $S$ to the simply-infinite set by $(\mathbb{N},0,S)$. 
Dedekind thus obtained the following axiomatic characterization of the natural numbers: \\
\ \\
{\bf Dedekind's Axioms} (Nr. 71 in Dedekind 1888 \cite{dedekind1888}) \\ 
The triple $(\mathbb{N},0,S)$ satisfies: \\
(D1) $0$ is not contained in the image of $S$. \\
(D2) $S$ is injective. \\
(D3) (Induction) The set $\mathbb{N}$ is simply-infinite, i.e.,
\[
\mathbb{N} = \bigcap_{C \, {\rm chain}, \,0 \in C} C = \{0,1:=S(0),2:=S(S(0)), \ldots \}.  
\]
\ \\
The following equivalent version is used by most authors \cite[\S 2]{sms} :  \\
\ \\
{\bf Dedekind-Peano Axioms} \\ 
The triple $(\mathbb{N},0,S)$ satisfies: \\
(DP1)  $0$ is not contained in the image of $S$. \\
(DP2)  $S$ is injective. \\
(DP3) (Induction) For any subset $M \subset \mathbb{N}$ such that
\begin{itemize}
\item[(i)] $0 \in M$. 
\item[(ii)] $n \in M \Rightarrow S(n) \in M$  
\end{itemize}
one has $M=\mathbb{N}$. \\
\ \\
Both descriptions above use quantification over subsets of $\mathbb{N}$, hence are formulated using a second-order logic framework.
Of course, Dedekind did not know about \emph{model theory} in the modern sense, 
and the distinction between syntax and semantics was not common then, but some of his remarks 
indicate that he took care to rule out possible \emph{non-standard models} with undesired properties
once some axioms would be relaxed. A posteriori, assuming a first-order logic framework, we know that (many) non-standard models indeed exist.

The most important result in Dedekind's book \cite{dedekind1888} is the famous \emph{recursion theorem}. It implies the 
\emph{categoricity} and \emph{semantical completeness} of $\mathbb{N}$. These two notions mean that all models of $\mathbb{N}$ are isomorphic
and any statement which holds for one model holds for any other model too. In modern form, the recursion theorem states: \\

{\bf Recursion Theorem} (Nr. 126 in Dedekind 1888 \cite{dedekind1888})\\
Let $\Omega$ be a set together with a self-map $\theta: \Omega \to \Omega$ und $\omega \in \Omega$. Then, there is one and only one map  
$\Psi: \mathbb{N} \to \Omega$, such that the diagram  
\begin{displaymath}
    \xymatrix{ \mathbb{N} \ar[d]_{S} \ar[r]^\Psi & \Omega \ar[d]^{\theta} \\
               \mathbb{N}  \ar[r]_{\Psi} & \Omega  }
\end{displaymath}
commutes, i.e., such that 
\[
\Psi(0)=\omega, \; \Psi(S(x))=\theta(\Psi(x)). 
\]
\ \\
There are more general versions, where the function $\theta$ depends on two variables and the rules are 
$\Psi(0)=\omega$ and $\Psi(S(x))=\theta(x,\Psi(x))$ \cite[\S 2]{sms}. 
Even more generally, we may look at recursion with additional auxiliary variables.
The recursion theorem in such a form was used in Dedekind's 1888 book \cite{dedekind1888} in \S 14 in the proof of Satz 159.  
It leads naturally to the definition of ``new'' functions $f={\rm rec}(g,h):\mathbb{N}^{n+1} \to \mathbb{N}$ 
by \emph{primitive recursion} \cite[\S 2]{sms}
\begin{eqnarray*}
f(0,y_1,\ldots,y_n) & =&g(y_1,\ldots,y_n) \\
f(S(x),y_1,\ldots,y_n) & =&h(x,f(x,y_1,\ldots,y_n),y_1,\ldots,y_n) 
\end{eqnarray*}
from given functions $g,h$.
General \emph{recursion theory} (alias the \emph{theory of computation}) needs one additional concept, called 
the (unbounded) search operator $\mu$, in order to obtain all \emph{computable} (alias \emph{partial recursive}) functions.  

The invention of recursion theory is often attributed to G\"odel or Skolem. However, it was Dedekind 
who first brought the recursive thinking to its full power by proving the recursion theorem.
Modern recursion theory owes to the work of Church, G\"odel, Herbrand, Hilbert (with Ackermann, Bernays), Grzegorczyk, 
K\'alm\'ar, Kleene, P\'eter, Post, Skolem, Turing et al. \cite[\S 6]{sms}.

\section{Reception} 

The reception of Dedekind's work is partly a story of misunderstandings, which had both positive and negative effects. 

Keferstein's critical, but mathematically not well-founded, comments on the 1888 book led to a famous letter of Dedekind \cite{sms,tapp}
which explains Dedekind's motivations and some ideas very well, including the potential existence of non-standard models of the natural numbers.  

The ``wrong'' proof of Satz 66 in Dedekind's 1888 book \cite{dedekind1888}, in which he claimed the existence of an infinite set 
by refering to non-mathematical objects (``thoughts''), was affected by the antinomies in set theory. Those eventually led to a crisis (``Grundlagenkrise'') 
in the foundations of mathematics after 1900. Dedekind also used the axiom of choice in his 1888 book. In contrast to the case of Frege, 
we do not know in detail how Dedekind thought and felt about the antinomies after his book was criticized on their basis. The only witness for this is  
Felix Bernstein who remembered that Dedekind told him that he had almost started to doubt the rationality of human thinking \cite[Vol. 3, p. 449]{dedekind}.
Nowadays, one postulates the existence of an infinite set and much of the foundational crisis seems to be exaggerated in retrospect. 
The antinomies had as an effect a lesser citation of Dedekind's work, for example by Hilbert and others, even though his work was silently incorporated in 
contemporary mathematics. In recent years, Dedekind's position in the foundations of mathematics is becoming readjusted again.

The beginnings of set theory were quite ``constructive'' with rather ``concrete'' sets.
There was no axiomatic method as in the later approach of Zermelo, who was a dedicated follower of Dedekind and called the axiom of infinity after him. 
Dedekind's letters with Cantor \cite{cavailles-noether} had a strong influence on the notions of cardinality and dimension in early set theory, reaching as far as to 
Brouwer's 1911 solution of the dimension problem for $\mathbb{R}^n$ as conjectured by Dedekind, i.e., that there is no homeomorphism between $\mathbb{R}^n$ 
and $\mathbb{R}^m$ unless $m=n$. 

There are also the letters to Lipschitz from 1876 \cite[Band 3, p. 468-479]{dedekind} in which Dedekind explained very well, and way before Hilbert and Peano, 
the essence of the axiomatic method. In addition, he made clear that the axiom of continuity is not part of Euklid's considerations, but rather a property of the ground field in question. 
Frege, Hilbert and Peano used the axiomatic method in a more consequent way, sometimes even without any reference to concrete models, e.g. in the case of 
Hilbert's axioms for the real numbers. In contrast, Dedekind required the set-theoretic existence at the beginning of any new creation (``Sch\"opfung'')
of a mathematical structure. 

Most surprising to me was Skolem's rejection of Dedekind's achievements in recursion theory \cite[\S 5]{sms}. In 1947, in a response to a remark 
in a paper of Curry, Skolem claimed to have discovered the recursive method in arithmetic in the year 1919 \cite[p. 499]{skolem}. Although Skolem had indeed invented the 
``primitiv recursive arithmetic'', a first-order fragment of Dedekind-Peano arithmetic avoiding quantifiers over infinite sets of integers, the recursive method 
of Dedekind and others was certainly still a major building block of this theory. 
Hao Wang was the first to correct this historical misinterpretation by stressing the importance of the earlier work of Dedekind, Peano and Gra{\ss}mann \cite[p. 17]{skolem}. 
As a pioneer for the introduction of first-order methods in mathematics, Skolem is, however, a very important figure.

\section{Influence in our times} 

Dedekind had a very strong influence on the development of algebraic number theory through the famous supplement \cite[Vol. 3]{dedekind} 
to Dirichlet's lectures on number theory, in which he developed algebraic number theory from scratch, using ideas of Galois and Kummer. 
His presentation became a raw model for most number theory lectures and books to follow. Hilbert, in his 
``Zahlbericht'', promoted Dedekind's theory even more. Subsequent developments include class field theory and later generalizations of the class number formula. 

There is a remarkable analogy between algebraic number theory and the theory of function fields of one variable. 
The latter was developed by Dedekind and H. Weber \cite[Vol. 1, p. 238]{dedekind}
in 1882 as well as independently by Kronecker in his ``Festschrift'' (1882) for Kummer. 
Dedekind wrote a small unpublished paper ``Bunte Bemerkungen'' on his readings of Kronecker's work.
From a modern perspective, both approaches are equivalent. 
The unification of algebraic number theory and the theory of function fields, which was intended by both parties, reaches out into the 20th century to 
Grothendieck's ``marriage'' of algebraic and arithmetic geometry via his theory of schemes and stacks. 
However, history has seen a lot of rivalry between followers of Dedekind and Kronecker. 
Weil mentions such a ``partisan war'' in his 1954 ICM address and strengthens the role of Kronecker. 

Emmy Noether adored Dedekind. This is apparent from her numerous comments in the 
collected works of Dedekind and her famous sentence ``Es steht alles schon bei Dedekind'', as witnessed by her student van der Waerden. 
This shows that she must have been on Dedekind's side in the partisan issue.

Dedekind's 1888 book \cite{dedekind1888} set the path for the axiomatic treatment of arithmetic, nowadays called \emph{Dedekind-Peano arithmetic}.
Recursion is also used in Hilbert's finitistic proof theory program as a metamathematical tool.
Furthermore, Dedekind's recursion theorem is related to G\"odel's two incompleteness theorems via the use of primitive recursive functions and 
(via transfinite or higher type versions) to the consistency proofs for arithmetic of Gentzen and G\"odel, hence to modern proof theory in general \cite[\S 6]{sms}.  
Finally, it affected of course the theory of computation, as well as intuitionistic and constructive mathematics. 
Algebraic set theory (Lawvere 1964) is a prototype example  for the influence of Dedekind on other mathematical foundations, including category theory and (higher) type theory.

Traces of Dedekind's presence can also still be felt in the famous \emph{Hilbert problems}, especially in 
Nr. 1 (continuum hypothesis), Nr. 2 (consistency of arithmetic), Nr. 7 (trancendence theory), 
Nr. 9 (reciprocity laws), Nr. 10 (diophantine equations), 
Nr. 11 (quadratic forms), Nr. 12 (generalized Kronecker-Weber theorem) and Nr. 24 (proof theory).
Among the Millenium problems of the Clay Foundation the problem ``P versus NP'' is very close to recursion theory 
and the Birch and Swinnerton-Dyer conjecture to algebraic number theory. 

Finally, the unpublished manuscript ``Sch\"opfung der Null und der negativen ganzen Zahlen'' (1872) may be viewed as being 
pointing to K-theory and the theory of motives of Grothendieck. 
Those theories are related to extensions of the class number formula of Dirichlet and Dedekind, e.g. Beilinson's conjectures and their variants, 
and are vital areas of research in modern arithmetic and
algebraic geometry. The theory of automorphic L-functions and the Langlands program can be seen as a vast
extension of Dedekind's $\zeta$-functions and the theory of modular forms.

\end{document}